\documentclass[journal,10pt]{IEEEtran}
\usepackage{graphicx,multirow,multicol,booktabs,pgfplots}

\usepackage{amsmath,bm}
\usepackage{amssymb,amsthm}
\usepackage[T1]{fontenc}
\usepackage[colorlinks,citecolor=blue]{hyperref}
\allowdisplaybreaks
\usepackage{setspace}
\usepackage{csvsimple}

\usepackage{amssymb}
\usepackage{amsmath}
\usepackage{url}
\usepackage{cite}

\usepackage{algorithm}
\usepackage[noend]{algpseudocode}

\algnewcommand\algorithmicinput{\textbf{Input:}}
\algnewcommand\Input{\item[\algorithmicinput]}

\algnewcommand\algorithmicoutput{\textbf{Output:}}
\algnewcommand\Output{\item[\algorithmicoutput]}

\usepackage[utf8]{inputenc}
\usepackage{color}

\usepackage{accents}

\usepackage{hyperref}

\allowdisplaybreaks
\theoremstyle{remark}

\pgfplotsset{compat=1.14}

\begin{document}
\title{\LARGE \bf Exact and Heuristic Approaches for the \\Stochastic $\bm N$-$\bm k$ Interdiction in Power Grids}

\author{Kaarthik Sundar$^{a}$, Andrew Mastin$^{b}$, Manuel Garcia$^{a}$, Russell Bent$^{a}$, Jean-Paul Watson$^{b}$
\thanks{$^{a}$ Los Alamos National Laboratory, Los Alamos, NM}
\thanks{$^{b}$ Lawrence Livermore National Laboratory, Livermore, CA}
\thanks{Corresponding author: \texttt{kaarthik@lanl.gov} }
\thanks{The authors acknowledge the funding provided by the U.S. Department of Energy's North American Energy Resilience Modeling (NAERM) program and the U.S. Department of Energy Office of Electricity's Advanced Grid Modeling program. The research work conducted at Los Alamos National Laboratory is done under the auspices of the National Nuclear Security Administration of the U.S. Department of Energy under Contract No. 89233218CNA000001. The work conducted at Lawrence Livermore National Laboratory is performed in part under the auspices of the U.S. Department of Energy under Contract DE-AC52-07NA27344.}
}

\markboth{}%
{}

\maketitle

\begin{abstract}
The article introduces the stochastic $\bm N$-$\bm k$ interdiction problem for power grid operations and planning that aims to identify a subset of $\bm k$ components (out of $\bm N$ components) that maximizes the expected damage, measured in terms of load shed. Uncertainty is modeled through a fixed set of outage scenarios, where each scenario represents a subset of components removed from the grid. We formulate the stochastic $\bm N$-$\bm k$ interdiction problem as a bi-level optimization problem and propose two algorithmic solutions. The first approach reformulates the bi-level stochastic optimization problem to a single level, mixed-integer linear program (MILP) by dualizing the inner problem and solving the resulting problem directly using a MILP solver to global optimality. The second is a heuristic cutting-plane approach, which is exact under certain assumptions. We compare these approaches in terms of computation time and solution quality using the IEEE-Reliability Test System and present avenues for future research. 
\end{abstract}

\IEEEpeerreviewmaketitle

\section{Introduction} \label{sec:intro}
As the linchpin of modern societies, power grids play a crucial role in sustaining socioeconomic systems. The escalating occurrence of natural disasters and deliberate attacks impacting power grids underscores the imperative to develop effective methods to pinpoint collections of critical components whose failure could lead to significant system-wide impacts. The deterministic $N$-$k$ interdiction problem \cite{salmeron2009worst} is one such method used to identify a small set of critical components in an electrical grid network that, when interdicted by an attacker, lead to a substantial amount of load shed. Here, $N$ denotes the total number of interdictable components in the system, and $k$ denotes the number of components that can be interdicted by the attacker. We interpret the problem as a Stackelberg game \cite{sundar2021credible} with one attacker (nature or an adversary) and one defender (system operator). The attacker first interdicts a set of $k$ components in the network; the defender responds by computing an operating point that minimizes the load shed. The deterministic $N$-$k$ interdiction problem seeks to determine an attack that maximizes this minimum load shed. Hence, the interdiction problem is naturally formulated as a bi-level optimization problem. This paper considers a stochastic version of the $N$-$k$ interdiction problem, where we assume that the on-off status of one or more of the existing components of the  deterministic $N$-$k$ interdiction problem is not known with certainty, and realizations of this uncertainty are made available through multiple scenarios. For instance, in a climate-driven extreme event setting, this uncertainty may correspond to the uncertainty associated with the evolution of the event and the components of the grid that could potentially impacted by the event. Each realization of the uncertainty (scenario) specifies a subset of components that are disabled, and we seek to determine an attack (set of $k$ components) that maximizes the expected minimum load shed, with the expectation computed as the average of the minimum load sheds across all scenarios. This results in the stochastic $N$-$k$ interdiction problem being formulated as a stochastic bi-level optimization problem. To the best of our knowledge, this is the first work in the literature to introduce, formulate and present solution approaches for the stochastic $N$-$k$ interdiction problem in power grids. 

The deterministic network interdiction problem has been well studied in the literature for both power grid \cite{salmeron2009worst,bienstock2010nk,sundar2021credible,mastin2023best} and more general capacitated flow networks \cite{brown2006defending,wood1993deterministic}. In contrast to capacitated flow networks, the power grid interdiction models include constraints representing the physics of power flow governed by non-convex alternating current (AC) power flow equations. For power grid interdiction problems, either convex relaxations of the AC power flow constraints \cite{sundar2018probabilistic} or linear direct current (DC) approximation of the same \cite{mastin2023best} are used to simplify the problem from both a computational and algorithmic standpoint. 

Although stochastic network interdiction has received substantial attention for capacitated flow networks \cite{cormican1998stochastic,held2005decomposition,janjarassuk2008reformulation}, it remains a largely unexplored topic for power grids. Interested readers are referred to a recent survey \cite{smith2020survey} for an exhaustive account of the different variants of the deterministic and stochastic network interdiction problems in the context of capacitated flow networks, with algorithmic approaches ranging from heuristic to exact methods. In this article, we first present the first bi-level stochastic optimization formulation for the stochastic $N$-$k$ interdiction problem with the DC power flow physics. The choice of DC power flow physics is motivated by two aspects: the first one being that the DC approximation is linear and the second is that it is often sufficient for transmission systems. Nevertheless, we remark that examining the impact of more accurate representation of power flow physics on the algorithms and the solutions is an important problem and we delegate this to future work. For the stochastic interdiction problem considering power grids, we propose two algorithmic approaches: (i) an exact approach based on reformulating the stochastic $N$-$k$ interdiction problem to a single-level MILP (extensive form MILP) through dualization of the inner problem and solving the extensive form MILP using an off-the-shelf commercial solver to global optimality, and (ii) an iterative heuristic alternating between solving the outer (attacker) and inner (defender) problems. The heuristic approach extends similar iterative algorithms previously introduced for the deterministic variant \cite{salmeron2009worst,sundar2021credible}. The heuristic approach is exact under certain assumptions \cite{sundar2021credible}. The motivation behind developing two approaches is as follows: The iterative algorithm is designed as a scalable heuristic for large-scale instances, whereas the exact algorithm has a global optimality guarantee and is used to quantify accuracy of the heuristic, but is only suitable in practice for smaller instances. We shall later test the validity of this hypothesis using extensive computational experiments. 

The rest of the article is organized as follows: in Sec. \ref{sec:formulation} we introduce notations and formulate the stochastic $N$-$k$ interdiction problem for power grids with the linear DC approximation of power flow physics. In Sec. \ref{sec:algo} and \ref{sec:results} we present an overview of the two algorithms we use to solve the problem and describe computational results that corroborate the effectiveness of these algorithms, respectively. We conclude by presenting avenues for future work in Sec. \ref{sec:conclusion}.

\section{Problem Formulation} \label{sec:formulation}
We now present necessary nomenclature and terminology to formulate the stochastic interdiction problem, including aspects that are well understood and used routinely
by the power systems community.   Without loss of generality, we assume that every bus has exactly one generator and one load.

\noindent \textit{Sets:} \\ 
$\mathcal B$ - buses (nodes) indexed by $b$ \\ 
$\mathcal L$ - transmission lines (edge) indexed by $(i, j)$\\ 
$\mathcal S$ - on-off scenarios indexed by $s$ \\ 
\noindent \textit{Variables:} \\ 
$x_{ij}$ - binary interdiction variable for line $(i, j)$ \\ 
$y_b^g$ - binary interdiction variable for the generator at bus $b$ \\
$\bar x$ - vector of interdiction variables $x_{ij}$ \\ 
$\bar y$ - vector of interdiction variables $y_b^g$ \\
$\eta_s(\bar x, \bar y)$ - auxiliary variable to model load shed given $s$, $\bar x$, $\bar y$ \\
$p_{ij}^s$ - power flow in line $(i, j)$ for scenario $s$ \\ 
$\theta_b^s$ - phase angle at bus $b$ for scenario $s$ \\
$\ell_b^s$ - percentage of load shed in bus $b$ for scenario $s$\\ 
$p_b^{gs}$ - generation at bus $b$ for scenario $s$ \\ 
\noindent \textit{Constants:} \\ 
$\bm b_{ij}$ - susceptance of line $(i, j)$ \\ 
$\bm t_{ij}$ - thermal limit of line $(i, j)$ \\ 
$\bm p_b^d$ - load on the bus $b$ \\ 
$\bm p_b^{gu}$ - maximum generation at bus $b$ \\
$\bm k$ - interdiction budget \\
$\bm \xi_{ij}^s$ - on/off (1/0) status of the line $(i, j)$ in scenario $s$ \\
$\bm \xi_b^{gs}$- on/off (1/0) status of the generator at bus $b$ in scenario $s$ \\
$\bm M$ - large constant value \\

The stochastic $N$-$k$ interdiction on a power grid can then be formulated as a bi-level stochastic optimization problem: 
\begin{subequations}
\begin{flalign}
& (\mathcal F_1) \quad z =  \max \, \frac 1{|\mathcal S|}\sum_{s \in \mathcal S} \eta_s(\bar x, \bar y) \text{ subject to:} & \label{eq:obj_outer} \\
& \sum_{(i,j) \in \mathcal E} x_{ij} + \sum_{b \in \mathcal N} y_b^g \leqslant \bm k, \text{ and } & \label{eq:outer-budget} \\ 
& (\bar x, \bar y) \in \{0, 1\}^{|\mathcal L| + |\mathcal B |} & \label{eq:outer-binary}
\end{flalign}
\label{eq:outer}
\end{subequations}
where 
\begin{subequations}
\begin{flalign}
& (\text{PLS}_s) \quad \eta_s(\bar x, \bar y) = \min \, \sum_{b \in {\mathcal B}} \bm p^d_b \ell_b^s \text{ subject to:}& \label{eq:obj_inner} \\
& p_b^{gs} - \bm p_b^d (1-\ell_b^s) = \sum_{i|(i,b) \in \mathcal L} p_{ib}^s - \sum_{i|(b, i) \in \mathcal L} p_{bi}^s ~ \forall b\in \mathcal B, & \label{eq:kcl} \\
& 0 \leqslant p^{gs}_b \leqslant (1-y_b^g)~\bm \xi_b^{gs}\, \bm p^{gu}_b  \quad \forall b \in \mathcal B, &\label{eq:pg} \\
& |p_{ij}^s|  \leqslant (1 - x_{ij})~\bm{\xi}_{ij}^s \bm t_{ij} \quad \forall (i,j) \in \mathcal L, & \label{eq:thermal} \\
& |p_{ij}^s + \bm b_{ij} (\theta^s_i - \theta^s_j)| \leqslant \bm M x_{ij} (1-\bm \xi_{ij}^s)  \quad \forall (i,j) \in \mathcal L, & \label{eq:pij} & \\
&0 \leqslant \ell_b^s \leqslant 1 \quad \forall b \in \mathcal B. & \label{eq:loadshed}
\end{flalign}
\label{eq:inner}
\end{subequations}
%

The outer problem in \eqref{eq:outer} is the attacker's problem and the inner problem in \eqref{eq:inner}, denoted by $\text{PLS}_s$ ($\text{PLS}$ abbreviates Primal Load Shedding problem), is the defender's problem. The inner problem \eqref{eq:inner} is specific to each scenario $s \in \mathcal S$; each inner problem decision variable is thus superscripted by $s$. The outer problem \eqref{eq:obj_outer} maximizes the expected objective of the inner problem across all scenarios, which we approximate via a sample average. The inner problem seeks to minimize the load shed given the interdicted components in scenario $s$. Constraints \eqref{eq:outer-binary} and \eqref{eq:outer-budget} enforce the binary restriction and budget constraint on the interdiction variables, respectively.  The constraints of the inner problem include (i) DC power flow constraints in \eqref{eq:kcl} and \eqref{eq:pij}, (ii) generation and line limits in \eqref{eq:pg} and \eqref{eq:thermal}, and (iii) bounds on load shed variables in \eqref{eq:loadshed}. We note that \eqref{eq:thermal} and \eqref{eq:pij} can be equivalently formulated as two linear constraints making the over all inner problem linear. 
Finally, while we include both standard bounds on line phase angle difference and bus-shunt models in our implementations, we omit these components in the inner problem \eqref{eq:inner} for brevity. 

\section{Algorithmic Approaches} \label{sec:algo}
We now present two approaches to solving the stochastic $N$-$\bm k$ interdiction problem for power grids. The first approach reformulates the bi-level stochastic optimization problem to a single-level stochastic optimization problem that is then solved using an off-the-shelf commercial MILP solver to global optimality. The second approach is a heuristic cutting-plane algorithm that is exact under certain assumptions. 

\subsection{Exact Approach} \label{subsec:exact}
We start out by formulating the dual problem for the inner optimization problem in $\text{PLS}_s$. To that end, we let $\pi_b^s$, $\varphi_b^{s}$, $\gamma_{ij}^{s+}$, $\gamma_{ij}^{s-}$, $\delta_{ij}^{s+}$, $\delta_{ij}^{s-}$, and $\omega_i^s$ denote the dual variables corresponding to \eqref{eq:kcl}, \eqref{eq:pg}, \eqref{eq:thermal}, \eqref{eq:pij}, and \eqref{eq:loadshed}, respectively. Then the dual of the inner problem $\text{PLS}_s$ is as follows: 
\begin{subequations}
\begin{flalign}
& (\text{DLS}_s) \quad \Delta_s(\bar x, \bar y) = \max \sum_{b \in \mathcal B} \left\{ (y_b^g-1)~\bm \xi_b^{gs}\, \bm p^{gu}_b \varphi_b^{s} -  \right.  & \notag \\
& \quad \left. \omega_b^s - \bm p_b^d \pi_b^s \right\} - \sum_{(i,j) \in \mathcal L} \left\{ \bm M x_{ij} (1-\bm \xi_{ij}^s) \left[ \delta_{ij}^{s+} + \delta_{ij}^{s-} \right] + \right. & \notag  \\
& \quad \left. \,(1 - x_{ij})~\bm{\xi}_{ij}^s \bm t_{ij} \left[ \gamma_{ij}^{s+} + \gamma_{ij}^{s-} \right] \right\} \text{ subject to: } & \label{eq:obj_dual} \\
& \pi_i^s - \pi_j^s + \gamma_{ij}^{s+} - \gamma_{ij}^{s-} + \delta_{ij}^{s+} - \delta_{ij}^{s-} = 0 ~~\forall (i,j) \in \mathcal L, & \label{eq:pij_d} \\
& \sum_{i|(i,b) \in \mathcal L} \bm b_{ib}\left[ \delta_{ib}^{s+} - \delta_{ib}^{s-} \right] - & \notag \\ 
& \qquad \qquad \sum_{i|(b,i) \in \mathcal L} \bm b_{bi}\left[ \delta_{bi}^{s+} - \delta_{bi}^{s-} \right] = 0 \quad \forall b \in \mathcal B, & \label{eq:theta_d} \\
& -\pi_b^s - \varphi_b^{s} = 0 \quad \forall b \in \mathcal B, & \label{eq:pg_d} \\
& -\bm p_b^d \pi_i^s - \omega_b^s \leqslant \bm p_b^d \quad \forall b \in \mathcal B, & \label{eq:l_d} \\
& \varphi_b^{s}, \omega_b^s \geqslant 0, ~~ \pi_b^s \text{ free } \quad \forall b \in \mathcal B, & \label{eq:nn_1} \\
& \gamma_{ij}^{s+}, \gamma_{ij}^{s-}, \delta_{ij}^{s+}, \delta_{ij}^{s-} \geqslant 0 \quad \forall (i,j) \in \mathcal L. & \label{eq:nn_2}
\end{flalign}
\label{eq:inner-dual}
\end{subequations}
Note that the solution space to the inner primal problem, $\text{PLS}_s$, is never empty for any interdiction plan $(\bar x, \bar y)$ and any scenario $s$ because, the solution where all generators do not generate any power is a feasible solution for $\text{PLS}_s$. We also remark that the objective value of the primal inner problem $\text{PLS}_s$ is bounded for all $(\bar x, \bar y) \in \{0, 1\}^{|\mathcal L| + |\mathcal B|}$. Since, $\text{PLS}_s$ is nonempty and bounded for all $(\bar x, \bar y)$ and $s$, the strong duality theorem of linear programming indicates that for any interdiction plan $(\bar x, \bar y)$ and scenario $s$, 
\begin{flalign}
    \eta_s(\bar x, \bar y) =  \Delta_s(\bar x, \bar y) \label{eq:strong-duality}
\end{flalign}
The above result yields the following single-level formulation for $\mathcal F_1$, which we refer to as $\mathcal F_2$.
\begin{subequations}
\begin{flalign}
& (\mathcal F_2) \quad z = \max \, \frac 1{|\mathcal S|}\sum_{s \in \mathcal S} \eta_s(\bar x, \bar y) \text{ subject to:}  & \label{eq:obj_outer_n} \\
& \eta_s(\bar x, \bar y) = \sum_{b \in \mathcal B} \left\{ (y_b^g-1)~\bm \xi_b^{gs}\, \bm p^{gu}_b \varphi_b^{s} - \omega_b^s - \bm p_b^d \pi_b^s \right\} -  & \notag \\
& \quad \sum_{(i,j) \in \mathcal L} \left\{ \bm M x_{ij} (1-\bm \xi_{ij}^s) \left[ \delta_{ij}^{s+} + \delta_{ij}^{s-} \right] + \right. & \notag  \\
& \quad \left. \,(1 - x_{ij})~\bm{\xi}_{ij}^s \bm t_{ij} \left[ \gamma_{ij}^{s+} + \gamma_{ij}^{s-} \right] \right\}, & \label{eq:sd_n} \\
& \text{Eq. \eqref{eq:outer-budget}, \eqref{eq:outer-binary}, and \eqref{eq:pij_d} -- \eqref{eq:nn_2} } \quad \forall s \in \mathcal S. & \notag 
\end{flalign}
\label{eq:single-nlp}
\end{subequations}
$\mathcal F_2$ is a mixed-integer nonlinear optimization problem. The nonlinearity in $\mathcal F_2$ arises due to binary-continuous product terms:  $y_b^g\varphi_b^{s}$, $x_{ij} \delta_{ij}^{s+}$, $x_{ij}\delta_{ij}^{s-}$, $x_{ij} \gamma_{ij}^{s+}$, and $x_{ij} \gamma_{ij}^{s-}$ in \eqref{eq:sd_n}; these bilinear terms can be equivalently reformulated into linear constraints using the well-known McCormick relaxation, which represents an exact reformulation for bilinear terms with a binary-continuous product \cite{mccormick1976computability}. Let $\mathcal F_3$ represent the resulting MILP obtained after applying Mccormick relaxation. The McCormick relaxation for $\sigma = B \cdot \mu$ where $B$ is a binary variable and $\mu \in [\mu^{\ell}, \mu^u]$ is a continuous variable is as follows:
\begin{subequations} \label{eq:McCormick}
    \begin{equation}
       B \mu^{\ell}  \leqslant \sigma \leqslant B \mu^{u} 
    \end{equation}
    \begin{equation}
      \mu + B \mu^u - \mu^u  \leqslant \sigma \leqslant \mu + B \mu^{\ell} - \mu^{\ell}
    \end{equation}
\end{subequations}
By introducing auxiliary variables 
\begin{flalign}
    \gamma_{ij}^s = \gamma_{ij}^{s+} + \gamma_{ij}^{s-} \text{ and } \delta_{ij}^s = \delta_{ij}^{s+} + \delta_{ij}^{s-} \label{eq:lifted-variables}
\end{flalign}
and applying the McCormick relaxation shown in \eqref{eq:McCormick} to the terms $y_b^g\varphi_b^{s}$, $x_{ij} \delta_{ij}^{s}$, and $x_{ij} \gamma_{ij}^{s}$, the formulation $\mathcal F_2$ can be reformulated into a single-level MILP. We remark that computing analytical lower and upper bounds on the dual variables $\varphi_b^{s}$, $\gamma_{ij}^{s}$ and $\delta_{ij}^s$ is an open question of interest and it greatly affects the convergence behaviour of off-the-shelf MILP solvers. For the purpose of this article, we assume provably valid bounds and delegate the question of computing tight bounds on these variables for future work. The exact approach presented is to directly use an MILP solver to solve $\mathcal F_3$ to global optimality. In the next section, we present a  brief overview of the heuristic approach to solve the stochastic $N$-$k$ interdiction problem.  

\subsection{Heuristic Cutting Plane Algorithm} \label{subsec:cp}
The heuristic is an iterative constraint generation algorithm that alternates between solving the outer and the inner problems and is designed so that convergence is efficiently. In the optimization literature, constraints that are linear and are iteratively added to the optimization problem are referred to as cutting planes and hence the terminology ``cutting-plane algorithm.'' The bi-level interdiction problem is solved using a constraint generation algorithm that alternates between solving the outer maximization problem in \eqref{eq:outer} and the inner problem $\text{PLS}_s$ for every scenario $s \in \mathcal S$. At each iteration, the algorithm relies on the construction of a linear function that upper-bounds the total load shed $\eta_s(\bar x, \bar y)$, given by the solution to the inner problem for any $s \in \mathcal S$. For any on-off scenario $s$ and an $N$-$k$ interdiction plan that satisfies the constraints in \eqref{eq:outer-budget} and \eqref{eq:outer-binary} denoted by $(\bar x^*, \bar y^*)$, let ${\eta}^*_s$ denote the minimum load shed that is obtained by solving the inner problem in \eqref{eq:inner} for the scenario $s$. Then, the algorithm uses the following constraint to upper-bound $\eta_s(\bar x, \bar y)$:
\begin{flalign}
    \eta_s(\bar x, \bar y) \leqslant \eta_s^* + \sum_{(i, j) \in \mathcal E} x_{ij} \alpha^s_{ij} + \sum_{b \in \mathcal B} y_b^g \beta^{gs}_b \label{eq:cut}
\end{flalign}
In each iteration, $|\mathcal S|$ number of such cuts, one for each on-off scenario, are added to the outer-problem in \eqref{eq:outer} and the outer problem is resolved to generate a new $N$-$k$ interdiction plan. The constraint \eqref{eq:cut} exists in the literature \cite{salmeron2009worst,sundar2021credible,sundar2018probabilistic} and found to be effective for the deterministic and probabilistic variants of the interdiction problem. The choice of coefficients $\alpha^s_{ij}$ and $\beta^{gs}_b$ presents a trade-off between the computation time and the quality of the interdiction plan provided by the algorithm. The equation used to compute these coefficients are as follows:
\begin{subequations}
    \begin{flalign}
        \alpha_{ij}^s &= \begin{cases}
            |p_{ij}^{s*}| & \text{if $\bm \xi_{ij}^s (1-x_{ij}^*)= 1$} \\ 
            0   & \text{otherwise}
        \end{cases} \label{eq:line-coeff} \\ 
        \beta^{gs}_b &= \begin{cases}
            p_{b}^{gs*} & \text{if $\bm \xi_{b}^{gs} (1 - y_b^{g*}) = 1$} \\ 
            0   & \text{otherwise}
        \end{cases} \label{eq:gen-coeff}
    \end{flalign}
    \label{eq:coeff}
\end{subequations}
Intuitively, using the coefficients \eqref{eq:coeff} in \eqref{eq:cut} has the following interpretation: if a line $(i, j) \in \mathcal E$ is
removed from a network, then the amount of load that is shed is at most the absolute value of power flowing on that line. When a generator is removed from the network, the amount of load shed is at most the amount of power that is being generated by that generator. In general, this is not always true due to Braess' paradox \cite{schafer2022understanding}. Hence, the constraint \eqref{eq:cut} that uses the coefficients in \eqref{eq:coeff} are not always valid for any interdiction plan and in that sense, the algorithm is a heuristic. Interested readers are referred to \cite{sundar2021credible} for techniques to compute coefficient values that are guaranteed to be valid; nevertheless as the authors in \cite{sundar2021credible} observe, using valid coefficients can lead to poor convergence of the algorithm and using coefficients in \eqref{eq:coeff} present a good trade-off between computation time and solution quality for solving the interdiction problem. Finally, for the sake of clarity, we present the pseudo-code of the heuristic in Algorithm \ref{algo:woods-algo}. In step \ref{step:cp} of the pseudo-code utilizing coefficients in \eqref{eq:coeff} may result in feasible solutions being removed by the outer problem solve in step \ref{step:outer-resolve} and $z_{\mathrm{ub}}$ in this case may be lower than the optimal objective value for the stochastic $N$-$k$ interdiction problem. 

\begin{algorithm}
\setstretch{1.1}
\caption{Pseudo-code for cutting-plane algorithm}\label{algo:woods-algo}
\begin{algorithmic}[1]
\vspace{1ex}
\Input tolerance, $\varepsilon > 0$
\Output $(\bar x^{\dagger}, \bar y^{\dagger})$
\State $z_{\mathrm{lb}} \gets -\infty$ \Comment{lower bound, feasible objective value}
\State $z_{\mathrm{ub}} \gets +\infty$ \Comment{upper bound} 
\State $(\bar x_I, \bar y_I)$ \Comment{initial $N$-$k$ attack}
\State \label{step:inner}solve inner problem $\forall s \in \mathcal S$ using $(\bar x_I, \bar y_I)$
\State update: $z_{\mathrm{lb}} \gets |\mathcal S|^{-1} \sum_{s \in \mathcal S} \eta_s(\bar x_I, \bar y_I)$ 
\State $\eta_s(\bar x_I, \bar y_I) \gets$ inner problem objective for each $s \in \mathcal S$
\State \label{step:cp}add constraint \eqref{eq:cut} for each $s \in \mathcal S$ to outer problem
\State \label{step:outer-resolve}resolve outer problem 
\State update: $z_{\mathrm{ub}} \gets \text{outer problem objective}$ 
\State update: $(\bar x_I, \bar y_I) \gets \text{outer problem solution}$ 
\If{$z_{\mathrm{ub}} - z_{\mathrm{lb}} \leqslant \varepsilon z_{\mathrm{lb}}$} output $(\bar x_I, \bar y_I)$ and stop  
\EndIf
\State return to step: \ref{step:inner}
\end{algorithmic}
\end{algorithm}

\section{Computational Results} \label{sec:results}
We now present extensive computational experiments to test the efficacy of both algorithmic approaches to solve the stochastic $N$-$k$ interdiction problem. We begin by providing the details of the data sets and the computational platforms used for the experiments. We use the IEEE Reliability Test System - Grid Modernization Lab Consortium (RTS-GMLC) with $73$ buses, $120$ transmission lines, and $96$ generators. The network data is open-sourced and can obtained from \url{https://github.com/GridMod/RTS-GMLC}. We remark that this test network is artificially geolocated and these locations are used to generate on-off scenarios for the interdiction problem. 

\begin{figure}
    \centering
    \includegraphics[scale=0.9]{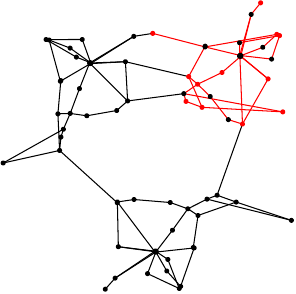}
    \caption{Buses and lines of the RTS-GMLC test system. For each scenario,  between $4$ to $6$ components (i.e., generators and lines) shown in red are randomly turned off. $200$ such scenarios are considered.}
    \label{fig:cluster}
\end{figure}

The on-off scenario sets $\mathcal S$ are generated as follows. First, a K-means clustering algorithm was run on the geo-located buses in the RTS-GMLC test system to obtain $3$ bus clusters. Then, one of the clusters is chosen and for each on-off scenario, $4$-$6$ components (generators in the buses and lines connecting the buses in the cluster) are randomly chosen based on a Bernoulli distribution to be turned off. The network, the generator buses, and the lines that are in the chosen cluster are shown in Fig. \ref{fig:cluster}  We generate $200$ such on-off scenarios and these scenarios constitute the set $\mathcal S$. The code for generating the scenarios and the scenarios themselves are open-sourced and made available through the GitHub repository: \url{https://github.com/kaarthiksundar/stochastic_nk}. 

Gurobi \cite{gurobi} was used to solve the linear and the MILP optimization problems for both algorithmic approaches. All experiments were run on an Intel Haswell 2.6 GHz, 62 GB, 20-core machine located at Los Alamos National Laboratory and finally, no time limit was imposed for any of the computational experiments. 

We first present a comparison of the computation times and the objective function values provided by the algorithmic approaches for attack budget $\bm k \in \{1, \dots, 10\}$. Table \ref{tab:times} reports both the objective value and the computation time in seconds for the heuristic and exact algorithms. It is clear from the table that the cutting-plane algorithm, in spite of having the ability to cut-off the global optimal solution, does not do so in practice and produces the same expected load shed (denoted by obj. in the table) values as that of the exact algorithm. As for the computation time, no algorithm between the two is a clear winner--the heuristic is able to converge to the global solution for $\bm k$ values $1-4$, $9$ and $10$ and the exact approach is observed to converge faster for the remaining values of $\bm k$. We also remark that one should not read too much into the computation times of the two approaches because they are not directly comparable in terms of the functionality because the exact approach is focused on proving optimality where as the heuristic approach is not. 

\begin{table}[!htb]
    \centering
        \small
    \caption{Comparison between the heuristic and exact algorithms in terms of solution quality and run-time in seconds. For all these runs, all the $200$ generated scenarios are used.}
    \label{tab:times}
    \begin{tabular}{crrrr}
        \toprule
        $\bm k$ & \multicolumn{2}{c}{Heuristic} & \multicolumn{2}{c}{Exact} \\ 
        \cmidrule{2-5}  
        & time (sec.) & obj. (MW) & time (sec.) & obj. (MW)\\ 
        \midrule
        \begin{filecontents*}{PMAPS-tables/times.csv}
k,cp_time,cp_obj,exact_time,exact_obj
1,60,151.82,860,151.82
2,86,481.32,1001,481.32
3,248,835.34,960,835.34
4,526,1190.34,1022,1190.34
5,2256,1545.34,1102,1545.34
6,6179,1900.34,2386,1900.34
7,4108,2255.34,1362,2255.34
8,4006,2605.34,1760,2605.34
9,1099,2942.59,2512,2942.59
10,236,3274.51,2800,3274.51
\end{filecontents*}
\csvreader[late after line=\\]{PMAPS-tables/times.csv}{1=\one,2=\two,3=\three,4=\four,5=\five}{\one & \two & \three & \four & \five}
        \bottomrule
    \end{tabular}
\end{table}

Table \ref{tab:vss} illustrates the value derived from solving the stochastic $N$-$k$ interdiction problem compared to an expected value problem (EVP) by computing the Value of Stochastic Solution (VSS) \cite{birge2011introduction}. In an EVP, all random variables are replaced by their expected values, and a deterministic $N$-$k$ interdiction is solved. The VSS is defined as the difference between the expected result of using the optimal EVP solution (denoted as EEV) and the objective of the stochastic $N$-$k$ interdiction problem, $z$. Generally, a higher VSS indicates greater value in using stochastic $N$-$k$ interdiction to model the problem instead of its deterministic counterpart. The last column in the table demonstrates that relative VSS values decrease as $\bm k$ increases. This observation aligns with the intuitive understanding that as a larger number of components, relative to the total number of components, are interdicted in the system, the value of uncertainty quantification diminishes.
\begin{table}[!htb]
    \centering
    \small
    \caption{Value of Stochastic Solution in (MW) for different $\bm k$ values.}
    \label{tab:vss}
    \begin{tabular}{rrrrr}
        \toprule
        $\bm k$ & $z$ (MW) & EEV (MW) & VSS (MW) & VSS (\%)\\ 
        \midrule
        \begin{filecontents*}{PMAPS-tables/VSS.csv}
k,RS,EEV,VSS,VSS_percent
1,151.82,30.5,121.32,79.91
2,481.32,458.47,22.85,4.75
3,835.34,633.21,202.13,24.2
4,1190.34,1033.21,157.13,13.2
5,1545.34,1388.21,157.13,10.17
6,1900.34,1839.41,60.93,3.21
7,2255.34,2232.26,23.08,1.02
8,2605.34,2569.51,35.83,1.38
9,2942.59,2904.41,38.18,1.3
10,3274.51,3241.66,32.85,1.0
\end{filecontents*}
\csvreader[late after line=\\]{PMAPS-tables/VSS.csv}{1=\one,2=\two,3=\three,4=\four,5=\five}{\one & \two & \three & \four & \five}
        \bottomrule
    \end{tabular}

\end{table}

Table \ref{tab:scenario_times} shows the increase in computation time for the heuristic approach for increasing sizes of $|\mathcal S|$. If computation time for large values of $|\mathcal S|$ is a concern, then Table \ref{tab:scenario_times} would help decide how to batch the scenarios in to smaller sets, solve each batch in parallel and obtain confidence intervals on the objective values and optimality gaps. This is a common approach used in the area of stochastic optimization to deal with the explosion of computation time with increasing cardinality of the scenario set \cite{janjarassuk2008reformulation}. 

\begin{table}[!htb]
    \centering
        \small
    \caption{Computation time in seconds taken by the heuristic for varying cardinality of the scenario set $|\mathcal S|$}
    \label{tab:scenario_times}
    \begin{tabular}{crrrr}
        \toprule
        $\bm k$ & $|\mathcal S| = 50$ & $|\mathcal S| = 100$ & $|\mathcal S| = 150$ & $|\mathcal S| = 200$ \\ 
        \midrule
        \begin{filecontents*}{PMAPS-tables/num_scenario_times.csv}
k,s50,s100,s150,s200
1,22,35,49,60
2,34,51,71,85
3,75,197,178,248
4,129,218,461,525
5,780,2018,3497,2256
6,981,2093,5212,6179
7,921,1580,3569,4107
8,625,1264,3465,4005
9,402,988,1904,1099
10,61,96,246,236
\end{filecontents*}
\csvreader[late after line=\\]{PMAPS-tables/num_scenario_times.csv}{1=\one,2=\two,3=\three,4=\four,5=\five}{\one & \two & \three & \four & \five}
        \bottomrule
    \end{tabular}
\end{table}

\begin{figure*}
    \centering
    \includegraphics[scale=0.93]{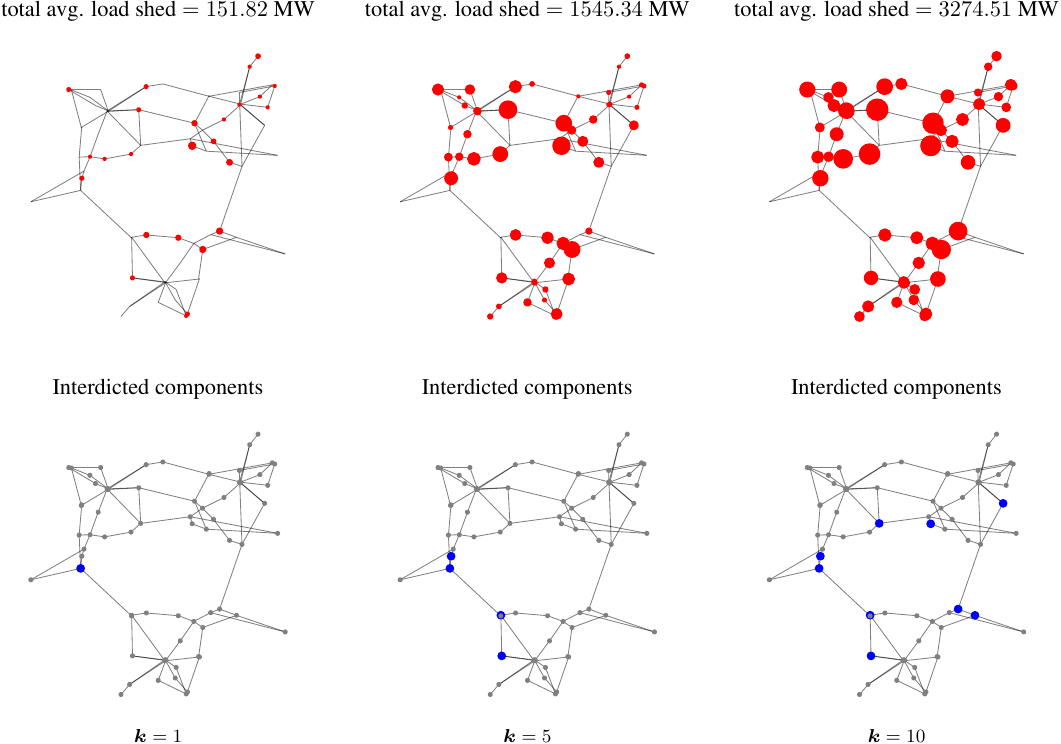}
    \caption{Top: Average load shed (over all the $200$ scenarios) at each bus in the optimal solution for the stochastic $N$-$k$ interdiction problem for different values of $\bm k$ is shown. Larger sizes of red discs indicate larger average load shed at a bus. Bottom: the components (shown in blue) are the interdicted components in the optimal $N$-$k$ interdiction plan.}
    \label{fig:solutions}
\end{figure*}

The subsequent set of results analyzes the tangible impact of the optimal $N$-$k$ interdiction plan on the network, quantified in terms of load loss at different locations. In this regard, Fig. \ref{fig:solutions} presents both the components in the optimal stochastic $N$-$k$ interdiction plan and the average load loss across the $200$ scenarios at each bus for $\bm k \in \{1, 5, 10\}$. In the figure, the size of the red disc at a bus indicates the magnitude of load loss – the larger the disc, the greater the load loss at that bus. The components in the optimal $N$-$k$ interdiction plan are highlighted in blue in the lower set of plots in Fig. \ref{fig:solutions}. Figures \ref{fig:cluster} and \ref{fig:solutions}, when considered together, provide a concrete depiction of how an operator would utilize the stochastic $N$-$k$ interdiction problem to identify vulnerabilities in the grid.

Assume that the components colored in red in Fig. \ref{fig:cluster} represent the set of components routinely affected by a specific weather-related extreme event, and that the actual components impacted by a particular event are uncertain. In such a scenario, the stochastic $N$-$k$ interdiction problem can be employed to identify the $\bm k$ most vulnerable components in the grid. These are the components whose removal during such an event would cause the maximum expected damage to the grid. The $\bm k$ most vulnerable components for $\bm k \in \{1, 5, 10\}$ and the corresponding maximum expected damage for the RTS-GMLC test network are illustrated in Fig. \ref{fig:solutions}.

\section{Conclusion and Avenues For Future Work} \label{sec:conclusion}
The stochastic interdiction problem in electric grids aims to pinpoint a limited set of vulnerable components in the grid, addressing the escalating challenges posed by extreme climate events and associated uncertainties. This article formulates the $N$-$k$ stochastic interdiction problem for electric transmission systems and introduces two algorithmic approaches: one heuristic and one exact. Computational experiments on the RTS-GMLC network demonstrate that the heuristic efficiently computes the globally optimal solution, despite lacking any global convergence guarantee. Numerous opportunities for future work are identified, including improving computation time performance for both approaches through (a) designing acceleration techniques and implementing decomposition algorithms like Progressive Hedging \cite{watson2011progressive}, commonly used in stochastic optimization, and (b) partitioning the larger scenario set into smaller batches to obtain confidence intervals on the objective values and optimality gaps. Exploring the scalability of the algorithms on practical power grids with tens of thousands of buses represents another promising avenue for future research. On the modeling front, particularly within a power systems context, it would be meaningful to introduce an additional layer of complexity into the existing interdiction problem by incorporating mitigation or protection measures for the defender.
\bibliographystyle{ieeetr}
\bibliography{references}

\end{document}